\documentclass[12pt]{article}
\newcommand{\be}{\begin{equation}}
\newcommand{\ee}{\end{equation}}
\newcommand{\bea}{\begin{eqnarray}}
\newcommand{\eea}{\end{eqnarray}}
\newcommand{\ba}{\begin{array}}
\newcommand{\ea}{\end{array}}

\newcommand{\bc}{\begin{center}}
\newcommand{\ec}{\end{center}}
\newcommand{\ben}{\begin{enumerate}}
\newcommand{\een}{\end{enumerate}}
\newcommand{\bfi}{\begin{figure}}
\newcommand{\efi}{\end{figure}}

\newcommand{\bq}{\begin{quote}}
\newcommand{\eq}{\end{quote}}
\newcommand{\bqu}{\begin{quotation}}
\newcommand{\equ}{\end{quotation}}
\newenvironment{emphit}{\begin{itemize}}{\end{itemize}}
\newcommand{\bemp}{\begin{emphit}}
\newcommand{\eemp}{\end{emphit}}

\newcommand{\bt}{\begin{tabular}}
\newcommand{\et}{\end{tabular}}

\newtheorem{myth}{Theorem}[section]
\newtheorem{mylem}{Lemma}[section]
\newtheorem{mycor}{Corollary}[section]

\newtheorem{mydef}{Definition}[section]
\newtheorem{myrem}{Remark}[section]
\newtheorem{myexam}{Example}[section]

\begin{document}
\date{}
\title{On A Pair Of Universal Weak Inverse
Property Loops \footnote{2000 Mathematics Subject Classification.
Primary 20NO5 ; Secondary 08A05}
\thanks{{\bf Keywords and Phrases :}weak inverse property loop, cross inverse property loop, ${\cal T}$
conditions, isotopism}}
\author{T\`em\'it\d{\'o}p\d{\'e} Gb\d{\'o}l\'ah\`an Ja\'iy\'e\d ol\'a\thanks{All correspondence to be addressed to this author.} \\
Department of Mathematics,\\
Obafemi Awolowo University, Ile Ife, Nigeria.\\
jaiyeolatemitope@yahoo.com,~tjayeola@oauife.edu.ng} \maketitle
\begin{abstract}
A new condition called ${\cal T}$ condition is introduced for the
first time and used to study a pair of isotopic loops. Under this
condition, a loop in the pair is a WIPL if and only if the other
loop is a WIPL. Furthermore, such WIPLs are isomorphic. The
translation elements $f$ and $g$ of a CIPL with the ${\cal T}$
condition(such that its $f,g$-isotope is an automorphic inverse
property loop) are found to be alternative, flexible, centrum and
equal elements. A necessary and sufficient condition for a pair
WIPLs with a weak ${\cal T}$ condition to be isomorphic is shown. A
CIPL and an isomorph are observed to have this weak ${\cal T}$
condition.
\end{abstract}

\section{Introduction}
\paragraph{}
Michael K. Kinyon \cite{phd33} gave a talk on Osborn Loops and
proposed the open problem : "Is every Osborn Loop universal?" which
is obviously true for universal WIP loops and universal CIP loops.
Our aim in this work is to introduced for the first time a new
condition called the ${\cal T}$ condition and use it to study a pair
of isotopic loops. The work is a special case of that of Osborn
\cite{phd89} for a WIPLs and we want to see if the result of Artzy
\cite{phd30} that isotopic CIP loops are isomorphic is true for WIP
loops or some specially related WIPLs(i.e special isotopes). The
special relation here is the ${\cal T}$ condition. But before these,
we shall take few basic definitions and concepts in loop theory
which are needed here. For more definitions, readers may check
\cite{phd3}, \cite{phd41},\cite{phd39}, \cite{phd49}, \cite{phd75}
and \cite{phd42}.

Let $L$ be a non-empty set. Define a binary operation ($\cdot $) on
$L$ : If $x\cdot y\in L$ for all $x, y\in L$, $(L, \cdot )$ is
called a groupoid. If the system of equations ;
\begin{displaymath}
a\cdot x=b\qquad\textrm{and}\qquad y\cdot a=b
\end{displaymath}
have unique solutions for $x$ and $y$ respectively, then $(L, \cdot
)$ is called a quasigroup. Furthermore, if there exists a unique
element $e\in L$ called the identity element such that for all $x\in
L$, $x\cdot e=e\cdot x=x$, $(L, \cdot )$ is called a loop. For each
$x\in L$, the elements $x^\rho ,x^\lambda\in L$ such that
$xx^\rho=e=x^\lambda x$ are called the right, left inverses of $x$
respectively. $L$ is called a weak inverse property loop (WIPL) if
and only if it obeys the weak inverse property (WIP);
\begin{displaymath}
xy\cdot z=e\Longrightarrow x\cdot yz=e~\forall~x,y,z\in L
\end{displaymath}
while $L$ is called a cross inverse property loop (CIPL) if and only
if it obeys the cross inverse property (CIP);
\begin{displaymath}
xy\cdot x^\rho=y.
\end{displaymath}
According to \cite{phd31}, the WIP is a generalization of the CIP.
The latter was introduced and studied by R. Artzy \cite{phd45} and
\cite{phd30} while the former was introduced by J. M. Osborn
\cite{phd89} who also investigated the isotopy invariance of the
WIP. Huthnance Jr. \cite{phd44} did so as well and proved that the
holomorph of a WIPL is a WIPL.  A loop property is called
universal(or at times a loop is said to be universal relative to a
particular property) if the loop has the property and every loop
isotope of such a loop possesses such a property. A universal WIPL
is called an Osborn loop in Huthnance Jr. \cite{phd44} but this is
different from the Osborn loop of Kinyon \cite{phd33} and Basarab.
The Osborn loops of Kinyon and Basarab were named generalised
Moufang loops or M-loops by Huthnance Jr. \cite{phd44} where he
investigated the structure of their holomorphs while Basarab
\cite{phd46} studied Osborn loops that are G-loops. Also,
generalised Moufang loops or M-loops of Huthnance Jr. are different
from those of Basarab. After Osborn's study of universal WIP loops,
Huthnance Jr. still considered them in his thesis and did an
elaborate study by comparing the similarities between properties of
Osborn loops(universal WIPL) and generalised Moufang loops. He was
able to draw conclusions that the latter class of loops is large
than the former class while in a WIPL the two are the same.

But in this present work, a new condition called ${\cal T}$
condition is introduced for the first time and used to study a pair
of isotopic loops. Under this condition, a loop in the pair is a
WIPL if and only if the other loop is a WIPL. Furthermore, such
WIPLs are isomorphic. The translation elements $f$ and $g$ of a CIPL
with the ${\cal T}$ condition(such that its $f,g$-isotope is an
automorphic inverse property loop) are found to be alternative,
flexible, centrum and equal elements. A necessary and sufficient
condition for a pair WIPLs with a weak ${\cal T}$ condition to be
isomorphic is shown. A CIPL and an isomorph are observed to have
this weak ${\cal T}$ condition.

\section{Preliminaries}
\begin{mydef}
Let $(L, \cdot )$ and $(G, \circ )$ be two distinct loops. The
triple $\alpha=(U, V, W) : (L, \cdot )\to (G, \circ )$ such that $U,
V, W : L\to G$ are bijections is called a loop isotopism
$\Leftrightarrow~ xU\circ yV=(x\cdot y)W~\forall ~x, y\in L$. Hence,
$L$ and $G$ are said to be isotopic whence, $G$ is an isotope of
$L$.

If $W=I$, then $\alpha $ is called a principal isotopism and in
addition, if $U=R_g$ and $V=L_f$ then $\alpha$ is called an
$f,g$-principal isotopism with the ordered pair $(g,f)$ called the
pair of translation elements of the principal isotope.
\end{mydef}

\begin{mydef}
Let $L$ be a loop. A mapping $\alpha\in S(L)$(where $S(L)$ is the
group of all bijections on $L$) which obeys the identity
$x^\rho=[(x\alpha)^\rho ]\alpha$ is called a weak right inverse
permutation. Their set is represented by $S_\rho (L)$.

Similarly, if $\alpha$ obeys the identity
$x^\lambda=[(x\alpha)^\lambda
 ]\alpha$ it is called a weak left inverse permutation. Their set is represented by $S_\rho (L)$

If $\alpha $ satisfies both, it is called a weak inverse
permutation. Their set is represented by $S'(L)$.

It can be shown that $\alpha\in S(L)$ is a weak right inverse if and
only if it is a weak left inverse permutation. So, $S'(L)=S_\rho
(L)=S_\lambda (L)$.
\end{mydef}

\begin{myrem}\label{0:1}
Every permutation of order 2 that preserves the right(left) inverse
of each element in a loop is a weak right (left) inverse
permutation.
\end{myrem}

\begin{myexam}\label{0:3}
If $L$ is an extra loop, the left and right inner mappings $L(x,y)$
and $R(x,y)~\forall~x,y\in L$ are automorphisms of orders 2
(\cite{phd36}). Hence, they are weak inverse permutations by
Remark~\ref{0:1}
\end{myexam}

\paragraph{}Throughout, we shall employ the use of the bijections;
$J_\rho~:~x\mapsto x^\rho$, $J_\lambda~:~x\mapsto x^\lambda$,
$L_x~:~y\mapsto xy$ and $R_x~:~y\mapsto yx$ for a loop and the
bijections; $J_\rho'~:~x\mapsto x^{\rho'}$, $J_\lambda'~:~x\mapsto
x^{\lambda'}$, $L_x'~:~y\mapsto xy$ and $R_x'~:~y\mapsto yx$ for its
loop isotope. If the identity element of a loop is $e$ then that of
the isotope shall be denoted by $e'$.

\begin{mylem}\label{0:2}
In a loop, the set of weak inverse permutations that commute form an
abelian group.
\end{mylem}

\begin{myrem}\label{0:4}
Applying Lemma~\ref{0:2} to extra loops and considering
Example~\ref{0:3}, it will be observed that in an extra loop $L$,
the Boolean groups $Inn_\lambda(L),Inn_\rho\le S'(L)$ .
$Inn_\lambda(L)$ and $Inn_\rho(L)$ are the left and right inner
mapping groups respectively. They have been investigated in
\cite{phd35}, and \cite{phd36}. This deductions can't be drawn for
CC-loops despite the fact that the left (right) inner mappings
commute and are automorphisms. And this is as a result of the fact
that the left(right) inner mappings are not of exponent 2.
\end{myrem}

\begin{mydef}(${\cal T}$-conditions)

Let $(G,\cdot )$ and $(H,\circ )$ be two distinct loops that are
isotopic under the triple $(A,B,C)$. $(G,\cdot )$ obeys the ${\cal
T}_1$ condition if and only if $A=B$. $(G,\cdot )$ obeys the ${\cal
T}_2$ condition if and only if
\begin{displaymath}
J_\rho'=C^{-1}J_\rho B=A^{-1}J_\rho C.
\end{displaymath}
$(G,\cdot )$ obeys the ${\cal T}_3$ condition if and only if
\begin{displaymath}
J_\lambda'=C^{-1}J_\lambda A=B^{-1}J_\lambda C.
\end{displaymath}
So, $(G,\cdot )$ obeys the ${\cal T}$ condition if and only if it
obey ${\cal T}_1$ and ${\cal T}_2$ conditions or ${\cal T}_1$ and
${\cal T}_3$ conditions since ${\cal T}_2\equiv {\cal T}_3$.
Furthermore, $(G,\cdot )$ obeys the ${\cal T}_{21}$ condition if and
only if $J_\rho'=C^{-1}J_\rho B$, $(G,\cdot )$ obeys the ${\cal
T}_{22}$ condition if and only if $J_\rho'=A^{-1}J_\rho C$,
$(G,\cdot )$ obeys the ${\cal T}_{31}$ condition if and only if
$J_\lambda'=C^{-1}J_\lambda A$ and $(G,\cdot )$ obeys the ${\cal
T}_{32}$ condition if and only if $J_\lambda'=B^{-1}J_\lambda C$.\\
So when $(H,\circ )=(G,\circ )$ is an $f,g$-principal isotope of
$(G,\cdot )$ under the triple $(R_g,L_f,I)$ : ${\cal T}_1$ condition
$\equiv R_g=L_f$, ${\cal T}_2$ condition $\equiv J_\rho'=J_\rho
L_f=R_g^{-1}J_\rho$, ${\cal T}_3$ condition $\equiv
J_\lambda'=J_\lambda R_g=L_f^{-1}J_\lambda$, ${\cal T}_{21}$
condition $\equiv J_\rho'=J_\rho L_f$, ${\cal T}_{22}$ condition
$\equiv J_\rho'=R_g^{-1}J_\rho$, ${\cal T}_{31}$ condition $\equiv
J_\lambda'=J_\lambda R_g$ and ${\cal T}_{32}$ condition $\equiv
J_\lambda'=L_f^{-1}J_\lambda$.

In case $(G,\cdot )$ and $(H,\circ )$ are two distinct non-isotopic
loops, then they are said to obey the weak ${\cal T}_{21}$ condition
if and only if $J_\rho'=A^{-1}J_\rho A$ or
$J_\lambda'=A^{-1}J_\lambda A$ for some $A~:~G\to H$ where
$J_\rho'$, $J_\rho$ and $J_\lambda'$, $J_\lambda$ still retain their
earlier definitions as right and left inverse mappings on $G$ and
$H$ respectively.
\end{mydef}
It must here by be noted that the ${\cal T}$-conditions refer to a
pair of isotopic loops at a time. This statement might be omitted at
times. That is whenever we say a loop $(G,\cdot )$ has the ${\cal
T}$-condition, then this is relative to some isotope $(H,\circ )$ of
$(G,\cdot )$

A loop $L$ is called a left inverse property loop(LIPL) if it obeys
the left inverse property (LIP);
\begin{displaymath}
x^\lambda (xy)=y~\forall~x,y\in L
\end{displaymath}
and a right inverse property loop(RIPL) if it obeys the right
inverse property (RIP);
\begin{displaymath}
(xy)y^\rho=x~\forall~x,y\in L.
\end{displaymath}
If it has both properties, then it is said to have the inverse
property (IP) hence called an inverse property loop (IPL).

\begin{mylem}\label{1:1}(\cite{phd3})
Let $L$ be a loop. The following are equivalent.
\begin{enumerate}
\item $L$ is a WIPL
\item $y(xy)^\rho =x^\rho~\forall~x,y\in L$.
\item $(xy)^\lambda x=y^\lambda~\forall~x,y\in L$.
\end{enumerate}
\end{mylem}

\begin{mylem}\label{1:2}
Let $L$ be a loop. The following are equivalent.
\begin{enumerate}
\item $L$ is a WIPL
\item $R_yJ_\rho L_y=J_\rho~\forall~y\in L$.
\item $L_xJ_\lambda R_x=J_\lambda~\forall~x\in L$.
\end{enumerate}
\end{mylem}

\section{Main Result}
\subsection{Isotopes of Weak Inverse Property Loops}
\begin{myth}\label{1:3}
Let $(G,\cdot )$ and $(H,\circ )$ be two distinct loops that are
isotopic under the triple $(A,B,C)$.
\begin{enumerate}
\item If the pair of $(G,\cdot )$ and $(H,\circ )$ obey the ${\cal T}$
condition, then $(G,\cdot )$ is a WIPL if and only if $(H,\circ )$
is a WIPL..
\item If  $(G,\cdot )$ and $(H,\circ )$ are WIPLs, then
\begin{displaymath}
J_\lambda R_xJ_\rho B=CJ_\lambda' R_{xA}'J_\rho'~\textrm{and}~J_\rho
L_xJ_\lambda A=CJ_\rho' L_{xB}'J_\lambda'~\forall~x\in G.
\end{displaymath}
\end{enumerate}
\end{myth}
{\bf Proof}\\
\begin{enumerate}
\item $(A, B,C)~:~G\rightarrow H$ is an isotopism $\Leftrightarrow xA\circ
yB=(x\cdot y)C\Leftrightarrow yBL_{xA}'=yL_xC\Leftrightarrow
BL_{xA}'=L_xC\Leftrightarrow L_{xA}'=B^{-1}L_xC\Leftrightarrow$
\begin{equation}\label{eq:1}
L_x=BL_{xA}'C^{-1}
\end{equation}
Also, $(A, B,C)~:~G\rightarrow H$ is an isotopism $\Leftrightarrow
xAR_{yB}'=xR_yC\Leftrightarrow AR_{yB}'=R_yC\Leftrightarrow
R_{yB}'=A^{-1}R_yC\Leftrightarrow$
\begin{equation}\label{eq:2}
R_y=AR_{yB}'C^{-1}
\end{equation}
Applying (\ref{eq:1}) and (\ref{eq:2}) to Lemma~\ref{1:2}
separately, we have : $R_yJ_\rho L_y=J_\rho $, $L_xJ_\lambda
R_x=J_\lambda \Rightarrow (AR_{xB}'C^{-1})J_\rho
(BL_{xA}'C^{-1})=J_\rho$, $(BL_{xA}'C^{-1})J_\lambda
(AR_{xB}'C^{-1})=J_\lambda\Leftrightarrow AR_{xB}'(C^{-1}J_\rho
B)L_{xA}'C^{-1}=J_\rho$, $BL_{xA}'(C^{-1}J_\lambda
A)R_{xB}'C^{-1}=J_\lambda\Leftrightarrow$
\begin{equation}\label{eq:3}
R_{xB}'(C^{-1}J_\rho B)L_{xA}'=A^{-1}J_\rho
C,~L_{xA}'(C^{-1}J_\lambda A)R_{xB}'=B^{-1}J_\lambda C.
\end{equation}

Let $J_\rho'=C^{-1}J_\rho B=A^{-1}J_\rho C$,
$J_\lambda'=C^{-1}J_\lambda A=B^{-1}J_\lambda C$. Then, from
(\ref{eq:3}) and by Lemma~\ref{1:2}, $H$ is a WIPL if $xB=xA$ and
$J_\rho'=C^{-1}J_\rho B=A^{-1}J_\rho C$ or $xA=xB$ and
$J_\lambda'=C^{-1}J_\lambda A=B^{-1}J_\lambda C \Leftrightarrow B=A$
and $J_\rho'=C^{-1}J_\rho B=A^{-1}J_\rho C$ or $A=B$ and
$J_\lambda'=C^{-1}J_\lambda A=B^{-1}J_\lambda C \Leftrightarrow A=B$
and $J_\rho'=C^{-1}J_\rho B=A^{-1}J_\rho C$ or
$J_\lambda'=C^{-1}J_\lambda A=B^{-1}J_\lambda C$. This completes the
proof of the forward part. To prove the converse, carry out the same
procedure, assuming the ${\cal T}$ condition and the fact that
$(H,\circ )$ is a WIPL.
\item If $(H,\circ )$ is a WIPL, then
\begin{equation}\label{eq:3.1}
R_y'J_\rho' L_y'=J_\rho',~\forall~y\in H
\end{equation}
while since $G$ is a WIPL,
\begin{equation}\label{eq:3.2}
R_xJ_\rho L_x=J_\rho~\forall~x\in G.
\end{equation}
The fact that $G$ and $H$ are isotopic implies that
\begin{equation}\label{eq:3.3}
L_x=BL_{xA}'C^{-1}~\forall~x\in G~and
\end{equation}
\begin{equation}\label{eq:3.4}
R_x=AR_{xB}'C^{-1}~\forall~x\in G.
\end{equation}
From (\ref{eq:3.1}),
\begin{equation}\label{eq:3.5}
R_y'=J_\rho' L_y'^{-1}J_\lambda'~\forall~y\in H~and
\end{equation}
\begin{equation}\label{eq:3.6}
L_y'=J_\lambda' R_y'^{-1}J_\rho'~\forall~y\in H
\end{equation}
while from (\ref{eq:3.2}),
\begin{equation}\label{eq:3.7}
R_x=J_\rho L_x^{-1}J_\lambda~\forall~x\in G~and
\end{equation}
\begin{equation}\label{eq:3.8}
L_x=J_\lambda R_x^{-1}J_\rho~\forall~x\in G.
\end{equation}
So, using (\ref{eq:3.6}) and (\ref{eq:3.8}) in (\ref{eq:3.3}) we get
\begin{equation}\label{eq:3.9}
J_\lambda R_xJ_\rho B=CJ_\lambda' R_{xA}'J_\rho'~\forall~x\in G
\end{equation}
while using (\ref{eq:3.5}) and (\ref{eq:3.7}) in (\ref{eq:3.4}) we
get
\begin{equation}\label{eq:3.10}
J_\rho L_xJ_\lambda A=CJ_\rho' L_{xB}'J_\lambda'~\forall~x\in G.
\end{equation}
\end{enumerate}

\begin{mycor}\label{1:5}
Let $(G,\cdot )$ and $(H,\circ )$ be two distinct loops that are
isotopic under the triple $(A,B,C)$. If $G$ is a WIPL with the
${\cal T}$ condition, then $H$ is a WIPL, :
\begin{enumerate}
\item there exists $\alpha ,\beta\in S'(G)$ i.e  $\alpha$ and $\beta $ are weak inverse
permutations and
\item $J_\rho'=J_\lambda'\Leftrightarrow J_\rho =J_\lambda$.
\end{enumerate}
\end{mycor}
{\bf Proof}\\
By Theorem~\ref{1:3}, $A=B$ and $J_\rho'=C^{-1}J_\rho B=A^{-1}J_\rho
C$ or $J_\lambda'=C^{-1}J_\lambda A=B^{-1}J_\lambda C$.
\begin{enumerate}
\item $C^{-1}J_\rho B=A^{-1}J_\rho C\Leftrightarrow J_\rho
B=CA^{-1}J_\rho C\Leftrightarrow J_\rho=CA^{-1}J_\rho CB^{-1}=
CA^{-1}J_\rho CA^{-1}=\alpha J_\rho \alpha $ where $\alpha
=CA^{-1}\in S(G,\cdot )$.
\item $C^{-1}J_\lambda A=B^{-1}J_\lambda C\Leftrightarrow J_\lambda
A=CB^{-1}J_\lambda C\Leftrightarrow J_\lambda = CB^{-1}J_\lambda
CA^{-1}=CB^{-1}J_\lambda CB^{-1}=\beta J_\lambda \beta$ where $\beta
=CB^{-1}\in S(G,\cdot )$.
\item $J_\rho'=C^{-1}J_\rho B$, $J_\lambda'=C^{-1}J_\lambda A$.
$J_\rho'=J_\lambda'\Leftrightarrow C^{-1}J_\rho B=C^{-1}J_\lambda
A=C^{-1}J_\lambda B\Leftrightarrow J_\lambda =J_\rho $.
\end{enumerate}

\begin{mylem}\label{1:6}
Let $(G,\cdot )$ be a WIPL with the ${\cal T}$ condition and
isotopic to another loop $(H,\circ )$. $(H,\circ)$ is a WIPL and $G$
has a weak inverse permutation.
\end{mylem}
{\bf Proof}\\
From the proof of Corollary~\ref{1:5}, $\alpha =\beta$, hence the
conclusion.

\begin{myth}\label{1:7}
With the ${\cal T}$ condition, isotopic WIP loops are isomorphic.
\end{myth}
{\bf Proof}\\
From Lemma~\ref{1:6}, $\alpha =I$ is a weak inverse permutation. In
the proof of Corollary~\ref{1:5}, $\alpha =CA^{-1}=I\Rightarrow
A=C$. Already, $A=B$, hence $(G,\cdot )\cong (H,\circ )$.

\subsection{$f,g$-Principal Isotopes of Weak Inverse Property Loops}

\begin{mylem}\label{1:3.4}
Let $(G,\cdot )$ be a WIPL with a WIP $f,g$-principal loop isotope
$(G,\circ )$ under the triple $\alpha =(R_g,L_f,I)$.
\begin{enumerate}
\item $J_\lambda R_fJ_\rho=L_f^{-1}$, $J_\rho
L_gJ_\lambda=R_g^{-1}$, $J_\lambda' R_g'J_\rho'=L_f$, $J_\rho'
L_f'J_\lambda'=R_g$.
\item $J_\rho L_f=R_f^{-1}J_\rho $, $J_\rho' L_f=R_g'J_\rho'$, $J_\lambda
R_g=L_g^{-1}J_\lambda$, $J_\lambda 'R_g=L_fJ_\lambda'$.
\item $J_\lambda R_f^{-1}J_\rho=J_\lambda' R_g'J_\rho'$, $J_\rho
L_g^{-1}J_\lambda=J_\rho' L_f'J_\lambda'$.
\item $\alpha =(J_\rho
L_g^{-1}J_\lambda ,J_\lambda R_f^{-1}J_\rho ,I)$ and $\alpha
=(J_\rho' L_f'J_\lambda',J_\lambda' R_g'J_\rho',I)$.
\end{enumerate}
\end{mylem}
{\bf Proof}\\
Using the second part of Theorem~\ref{1:3},
\begin{equation}\label{eq:11}
J_\lambda R_xJ_\rho L_f=J_\lambda' R_{xg}'J_\rho'~and
\end{equation}
\begin{equation}\label{eq:12}
J_\rho L_xJ_\lambda R_g=J_\rho' L_{fx}'J_\lambda'~\forall~x\in G.
\end{equation}
(1) to (4) are achieved by taken $x=f$, $x=e$ in (\ref{eq:11}) and
$x=g$, $x=e$ in (\ref{eq:12}).

\begin{myrem}
In [I.4.1~Theorem, \cite{phd3}], it is shown that in an IP
quasigroup ; $J_\lambda R_xJ_\rho=L_{x^\lambda}$, $J_\rho
L_xJ_\lambda=R_{x^\rho}$. These equations are true relative to the
translational elements $f,g$ in a universal WIPL as shown in
Lemma~\ref{1:3.4} (1).
\end{myrem}

\begin{mycor}\label{1:3.5}
Let $(G,\cdot )$ be a WIPL with an $f,g$-principal loop isotope
$(G,\circ )$. If the ${\cal T}$ condition holds in $(G,\cdot )$,
then $(G,\circ )$ is a WIPL. But, provided any of the following
holds :
\begin{enumerate}
\item ${\cal T}_1$ and ${\cal T}_{21}$ conditions
\item ${\cal T}_1$ and ${\cal T}_{22}$ conditions,
\end{enumerate}
if $(G,\circ )$ is a WIPL, then the ${\cal T}$ condition holds.
Hence, $L_f,R_g\in S'(G,\cdot )$.
\end{mycor}
{\bf Proof}\\
The proof of the first part is like the proof of the first part of
Theorem~\ref{1:3}. The second part is achieved using (2) of
Lemma~\ref{1:3.4}.

\begin{mycor}\label{1:4}
Let $(G,\cdot )$ be a loop with an $f,g$-principal isotope $(G,\circ
)$. If~ $(G,\cdot )$ is a WIPL, then $(G,\circ )$ is a WIPL provided
the ${\cal T}$ condition holds in $(G,\cdot )$ . But, if $(G,\circ
)$ is a WIPL then $J_\lambda R_xJ_\rho L_f=J_\lambda'
R_{xg}'J_\rho'$ and $J_\rho L_xJ_\lambda R_g=J_\rho'
L_{fx}'J_\lambda'~\forall~x\in G$.
\end{mycor}
{\bf Proof}\\
In Theorem~\ref{1:3}, let $(A,B,C)=(R_g,L_f,I)$ and $G=H$, then
$A=R_g,B=L_f$ and $C=I$. Putting these in the results ;
$A=R_g=B=L_f$ and $J_\rho'=I^{-1}J_\rho L_f=R_g^{-1}J_\rho I$ or
$J_\lambda'=I^{-1}J_\lambda R_g=L_f^{-1}J_\lambda I\Leftrightarrow
R_g=L_f$ and $J_\rho'=J_\rho L_f=R_g^{-1}J_\rho $ or
$J_\lambda'=J_\lambda R_g=L_f^{-1}J_\lambda$. This completes the
proof of the first part. The second part follows by just using the
above replacements.

\begin{mycor}\label{1:4.1}
Let $G$ be a WIPL with either the ${\cal T}_2$ or ${\cal T}_3$
condition. For any arbitrary $f,g$-principal isotope $G'$ of $G$,
the principal isotopism is described by the triple $(J_\rho
J_\lambda',J_\lambda J_\rho',I)$.
\end{mycor}
{\bf Proof}\\
By the ${\cal T}_2$ or ${\cal T}_3$ condition, $J_\rho'=J_\rho L_f$,
$J_\lambda'=J_\lambda R_g\Rightarrow L_f=J_\lambda J_\rho'$,
$R_g=J_\rho J_\lambda'$. Thus the triple $(R_g,L_f,I)=(J_\rho
J_\lambda',J_\lambda J_\rho',I)$.

\begin{mylem}\label{1:4.2}
For every WIPL with the ${\cal T}_2$ condition, if all
$f,g$-principal isotopes have the same right(left) inverse mappings,
then there exists a unique $f,g$-principal loop isotope. Thence, all
loop isotopes of such a loop are isomorphic loops. Furthermore, with
the ${\cal T}$ condition, if all $f,g$-principal isotopes have the
same right(left) inverse mappings then, there exists a unique W. I.
P. $f,g$-principal loop isotope. Thence, all loop isotopes of such a
loop are isomorphic WIP loops.
\end{mylem}
{\bf Proof}\\
Let $G$ be the WIPL in consideration. Let $G'$ and $G''$ be any two
distinct principal isotopes of $G$ under the triples $\alpha
=(R_{g'},L_{f'},I)$ and $\beta =(R_{g''},L_{f''},I)$ respectively.
If the ${\cal T}_2$ condition holds in $G$, then $L_{f'}=J_\lambda
J_\rho'$, $R_{g'}=J_\rho J_\lambda'$ and $L_{f''}=J_\lambda
J_\rho''$, $R_{g''}=J_\rho J_\lambda''$. By hypothesis, $J_\rho
=J_\rho''$ or $J_\lambda =J_\lambda''$. So, $R_{g'}=R_{g''}$ and
$L_{f'}=L_{f''}\Rightarrow f'=f''$ and $g'=g''$. Thus, $\alpha
=\beta\Rightarrow G'=G''$. This proves the uniqueness of
$f,g$-principal loop isotope. Recall that if $H'$ and $H''$ are any
two distinct loop isotopes of $G$ then there exists $G'$ and $G''$
such that $H'\cong G'$ and $H''\cong G''$. So, $H'\cong H''$.

Assuming the ${\cal T}$ condition in $G$, the further statement
follows by the same argument and the isotopes are therefore WIP
loops since this property is isomorphic invariant.

\paragraph{Centrum}The set of elements that commute with all other elements in a
loop $L$ is denoted by $C(L)$.

\begin{mylem}\label{1:8}
\begin{enumerate}
\item Let $(G,\cdot )$ be a WIPL with the ${\cal T}$ condition for an $f,g$-principal isotope $(G,\circ
)$ or
\item if an $f,g$-principal isotope of a WIPL $(G,\cdot )$ with the ${\cal T}_1$ and ${\cal T}_{21}$
conditions or ${\cal T}_1$ and ${\cal T}_{22}$ conditions is also a
WIPL, then
\end{enumerate}
\begin{description}
\item[(a)] $xg=fx$, $f,g\in C(G)$.
\item[(b)] $x^{\rho'}=fx^\rho$.
\item[(c)] $x^{\lambda'}=x^\lambda g$.
\item[(d)] $gg=ff=fg=gf$.
\item[(e)] $f^{\rho'}=g^{\lambda'}=e$
\end{description}
\end{mylem}
{\bf Proof}\\
Using Corollary~\ref{1:3.5} ;
\begin{description}
\item[(a)] $R_g=L_f\Leftrightarrow xR_g=xL_f\Leftrightarrow
xg=fx$. \item[(b)] $x^{\rho'}=xJ_\rho'=xJ_\rho L_f=x^\rho
L_f=fx^\rho$. \item[(c)] $x^{\lambda'}=xJ_\lambda'=xJ_\lambda
R_g=x^\lambda R_g=x^\lambda g$. \item[(d)] From (a), with $x=f$,
$fg=ff$. With $x=g$, $gg=fg\Rightarrow fg=ff=gg$. \item[(e)] Putting
$x=f$ in (b), $f^{\rho'}=ff^\rho =e$. Putting $x=g$ in (c),
$g^{\lambda'}=g^\lambda g=e$.
\end{description}

\begin{mycor}\label{1:8.5}
If an $f,g$-principal isotope of a CIPL $G$ with the ${\cal T}$
condition is an AIPL, then $f$ and $g$ are
\begin{enumerate}
\item Alternative elements (i.e $(xx)y=x(xy)$ and
$y(xx)=(yx)x~\forall~y\in G$ and $x\in \{f,g\}$).
\item Flexible elements (i.e $x(yx)=(xy)x~\forall~y\in G$ and $x\in
\{f,g\}$).
\item Centrum elements (i.e $xy=yx~\forall~y\in G$ and $x\in
\{f,g\}$).
\item Equal elements (i.e $f=g$).
\end{enumerate}
\end{mycor}
{\bf Proof}\\
These are achieved using Lemma~\ref{1:8} and some results in
\cite{phd30}.

\begin{myrem}
The properties of $f$ and $g$ proved in Corollary~\ref{1:8.5} above
were not gotten in \cite{phd30}. In fact, in addition to some
identities and autotopisms stated in \cite{phd30} satisfied by $f$
and $g$, we have the following :

\paragraph{Identities;}
$$xg\cdot gy=gg\cdot xy=(xy\cdot g)g=(g\cdot xy)g=ff\cdot
xy=f(xy\cdot g)=f(f\cdot xy).$$
$$xg\cdot g=g\cdot gx=gg\cdot x=xg\cdot g=gx\cdot g=ff\cdot x=f\cdot
xg=f\cdot fx.$$ for all $x$ and $y$ in a CIPL with the ${\cal T}$
condition.

\paragraph{Autotopisms;}
$$(R_g,L_g,L_{gg}),(R_g,L_g,R_g^2),(R_g,L_g,L_gR_g),(L_f,L_g,L_fR_g)$$
$$(L_f,L_g,L_gL_f), (L_f,L_g,R_gL_f),(L_f,L_g,L_f^2)(L_f,L_g,L_{ff})$$
for translation elements $f$ and $g$.
\end{myrem}

\begin{mycor}\label{1:9}
If a WIPL $G$ has the ${\cal T}$ condition, then
\begin{description}
\item[(a)] $f^{\rho'}=g^{\lambda'}$ in $G'$.
\item[(b)] $gg=ff=e'$ in $G'$.
\end{description}
\end{mycor}
{\bf Proof}\\
These follows from Lemma~\ref{1:8}.
\begin{description}
\item[(a)] This follows immediately by (e);
$f^{\rho'}=e=g^{\lambda'}$. \item[(b)] In $G'$, $fg$ is the identity
element. Thus from (d), $gg=ff=fg=e'$.
\end{description}

\paragraph{Automorphic Inverse Property Elements}
In a loop $L$, an element $g\in L$ is said to be a right automorphic
inverse property element ($\rho$-AIPE) if and only if
$(xg)^{-1}=x^{-1}g^{-1}~\forall~x\in L$ while $g$ is called a left
automorphic inverse property element ($\lambda$-AIPE) if and only if
$(gx)^{-1}=g^{-1}x^{-1}$ hence $g$ is called an automorphic inverse
property element if and only if these two conditions hold together.

We also use these definitions for anti-automorphic inverse property
elements(AAIPE) by simply interchanging the positions of $g^{-1}$
and $x^{-1}$ on the right hand side.

\begin{mylem}\label{1:11}
A LIPL(RIPL) is a WIPL if and only if it is a RIPL(LIPL).
\end{mylem}

\begin{myth}\label{1:12}
If an $f,g$-principal isotope of a LIP(RIP), WIPL $G$ is again a
WIPL, then
\begin{description}
\item[(a)] $g$ is a $\rho$-AIPE
\item[(b)] $f$ is a $\lambda$-AIPE
\item[(c)] $g,f\in C(G)$.
\end{description}
\end{myth}
{\bf Proof}\\
As shown in Lemma~\ref{1:11}, a LIPL(RIPL) is a W. IPL if and only
if it is a RIPL(LIPL). Using Corollary~\ref{1:4} and Lemma~\ref{1:8}
;
\begin{description}
\item[(a)] $x^{\rho'}=fx^\rho =(xg^{-1})^\rho\Rightarrow
fx^{-1}=(xg^{-1})^{-1}$. So,
$fx^{-1}=x^{-1}g=(xg^{-1})^{-1}\Rightarrow x^{-1}g^{-1}=(xg)^{-1}$.
\item[(b)] $x^{\lambda'}=x^\lambda g=(f^{-1}x)^\lambda \Rightarrow
x^{-1}g=(f^{-1}x)^{-1}$. So,
$x^{-1}g=fx^{-1}=(f^{-1}x)^{-1}\Rightarrow f^{-1}x^{-1}=(fx)^{-1}$.
\item[(c)] By (a), $x^{-1}g^{-1}=(xg)^{-1}=g^{-1}x^{-1}\Rightarrow x^{-1}g^{-1}=g^{-1}x^{-1}\Rightarrow
xg=gx\Rightarrow g\in C(G)$.

Similarly by (b), $f^{-1}x^{-1}=(fx)^{-1}=x^{-1}f^{-1}\Rightarrow
f^{-1}x^{-1}=x^{-1}f^{-1}\Rightarrow xf=fx\Rightarrow f\in C(G)$.
\end{description}

\begin{myrem}
In \cite{phd12}, it is shown that if a loop $L$ has the RIP(L. I.
P.) then,
\begin{description}
\item[(i)] $(xy)^{-1}=y^{-1}x^{-1}~\forall~y\in N_\rho
(L)\big((N_\lambda (L)\big)$ and $x\in L$
\item[(ii)] $(yx)^{-1}=x^{-1}y^{-1}~\forall~y\in N_\rho
(L)\big((N_\lambda (L)\big)$ and $x\in L$
\end{description}
where $N_\rho (L)\big((N_\lambda (L)\big)$ is the right(left)
nucleus of $L$. This implies right(left) nuclear elements are both
$\rho$-AAIPE and $\lambda$-AAIPE.
\end{myrem}

\subsection{Isomorphic Weak Inverse Property Loops}
\begin{myth}\label{1:14}
Let $(G,\cdot )$ and $(H,\circ )$ be two distinct WIP loops with the
weak ${\cal T}_{21}$ condition. $(G,\cdot )\cong (H,\circ )$  if and
only if $AJ_\rho B=CJ_\rho D\Rightarrow A=C$ or $B=D$ for some
$A,B,C,D\in S(G,\cdot )$ where $J_\rho$ is the right inverse mapping
on $(G,\cdot )$.
\end{myth}
{\bf Proof}\\
$G$ is a WIPL $\Leftrightarrow $
\begin{equation}\label{eq:4}
R_xJ_\rho L_x=J_\rho~\forall~x\in G~while
\end{equation}
$H$ is a W. I. P L. $\Leftrightarrow R_y'J_\rho'
L_y'=J_\rho'~\forall~y\in H$. Let $y=xA$ such that $A~:~(G,\cdot
)\rightarrow (H,\circ )$ is a bijection. Let $J_\rho'=A^{-1}J_\rho
A$, then  $R_y'(A^{-1}J_\rho A)L_y'=A^{-1}J_\rho A\Leftrightarrow
AR_y'(A^{-1}J_\rho A)L_y'A^{-1}=J_\rho\Leftrightarrow $
\begin{equation}\label{eq:5}
(AR_y'A^{-1})J_\rho (AL_y'A^{-1})=J_\rho
\end{equation}
Using the hypothesis combined with (\ref{eq:4}) and (\ref{eq:5}) ;
$R_x=AR_{xA}'A^{-1}$ or $L_x=AL_{xA}'A^{-1}\Leftrightarrow (G,\cdot
)\cong (H,\circ )$.

Conversely ; $G$ is a WIPL $\Leftrightarrow $
\begin{equation}\label{eq:6}
R_xJ_\rho L_x=J_\rho~\forall~x\in G
\end{equation}
$H$ is a WIPL $\Leftrightarrow $
\begin{equation}\label{eq:7}
R_y'J_\rho' L_y'=J_\rho'~\forall~y\in H
\end{equation}
If $G\cong H$, then $\exists $ a bijection $A~:~G\rightarrow
H~\ni~xA\circ yA=(x\cdot y)A~\forall~x,y\in G\Leftrightarrow $
\begin{equation}\label{eq:8}
L_x=AL_{xA}'A^{-1}~or~R_y=AR_{yA}'A^{-1}~\forall~x,y\in G
\end{equation}
Putting $J_\rho'=A^{-1}J_\rho A$ in (\ref{eq:7}), we have
$R_y'(A^{-1}J_\rho A)L_y'=A^{-1}J_\rho A\Leftrightarrow
(AR_y'A^{-1})J_\rho (AL_y'A^{-1})=J_\rho$. Let
$B=AR_y'A^{-1},C=AL_y'A^{-1}\in S(G,\cdot )$ and $D=R_x,E=L_x\in
S(G,\cdot )$, then
\begin{equation}\label{eq:9}
BJ_\rho C=J_\rho
\end{equation}
\begin{equation}\label{eq:10}
DJ_\rho E=J_\rho
\end{equation}
Thus, (\ref{eq:9}) and (\ref{eq:10}) implies (\ref{eq:8}) i.e
$BJ_\rho C=J_\rho$ and $DJ_\rho E=J_\rho\Rightarrow B=D$ or $C=E$.
Hence the proof.

\begin{myrem}
By Theorem~\ref{1:14}, WIP loops with the weak ${\cal T}_{21}$
condition are isomorphic under a necessary and sufficient condition.
This condition is therefore an isomorphy condition and not an
isotopy-isomophy condition which is characteristic of the Wilson's
Identity(\cite{phd90}) and the condition given in [Lemma~2,
\cite{phd89}] by Osborn.
\end{myrem}

\begin{myth}\label{2:13}
CIP is an isomorphic invariant property if and only if
$CJ_\rho'=J_\rho C$ or $DJ_\lambda'=J_\lambda D$ where $C$ and $D$
are permutations while $J_\rho, J_\lambda$ and $J_\rho', J_\lambda'$
are the right, left inverse mappings of the loop and its isomorph
respectively.
\end{myth}
{\bf Proof}\\
Let $(G,\cdot )$ be a CIPL isomorphic to a loop $(H,\circ )$ i.e
there exists $A~:~G\to H~\ni~G\cong^{^A}H$.
\begin{equation}\label{eq:13}
xA\circ yA=(x\cdot y)A\Leftrightarrow R_{yA}'=A^{-1}R_yA
\end{equation}
\begin{equation}\label{eq:14}
xA\circ yA=(x\cdot y)A\Leftrightarrow L_{xA}'=A^{-1}L_xA
\end{equation}
$G$ is a CIPL $\Leftrightarrow xy\cdot x^\rho\Leftrightarrow
R_{x^\rho }=L_x^{-1}$. Let $z=x^\rho$, then $x=z^\lambda$. So,
$R_z=L_{z^\lambda }^{-1}$. Putting this in (\ref{eq:13}) ;
\begin{equation}\label{eq:15}
R_{yA}'=A^{-1}L_{y^\lambda }^{-1}A
\end{equation}
From (\ref{eq:14}), $L_{x^\lambda }=AL_{x^\lambda A}'A^{-1}$.
Putting this in (\ref{eq:15}), get $R_{yA}'=A^{-1}L_{y^\lambda
}^{-1}A=A^{-1}(AL_{x^\lambda A}'A^{-1})^{-1}A=L_{y^\lambda
}A'^{-1}\Rightarrow L_{y^\lambda }A'R_{yA}'=I\Leftrightarrow
zAL_{y^\lambda }A'R_{yA}'=zA\Leftrightarrow (yJ_\lambda A\circ
zA)\circ yA=zA$. Let $C=J_\lambda A\Leftrightarrow A=J_\rho C$, so
\begin{equation}\label{eq:16}
(yC\circ zA)\circ yJ_\rho C=zA
\end{equation}
If CIP is an isomorphic invariant property then $(y'\circ z')\circ
y'^{\rho'}=z'~\forall~y',z'\in H$. Thus, with $y'=yC$, $z'=zA$,
\begin{equation}\label{eq:17}
(yC\circ zA)\circ yCJ_\rho'=zA
\end{equation}
Comparing (\ref{eq:16}) and (\ref{eq:17}), $J_\rho C=CJ_\rho'$.

The converse is proved by doing the reverse. The proof of the second
part is similar.

End of proof.

\begin{mycor}\label{2:13.1}
A CIPL and a loop isomorph have the weak ${\cal T}_{21}$ condition.
\end{mycor}
{\bf Proof}\\
Since the CIP is an isomorphic invariant property, the proof of this
claim follows from Theorem~\ref{2:13}.

\begin{mylem}\label{2:14}
Let $G$ be a CIPL. Then, $D=J_\rho^2C$ and $C=J_\lambda^2D$. If in
addition the loop has the RIP or LIP, then $D=C$. Hence, $J_\rho
=J_\lambda$.
\end{mylem}
{\bf Proof}\\
Using the notations in Theorem~\ref{2:13}, since $C=J_\lambda A$ and
$D=J_\rho A$ implies that $A=J_\lambda^{-1}C=J_\rho C$ or
$A=J_\rho^{-1}D=J_\lambda D$ then $D=J_\rho^2C$ and
$C=J_\lambda^2D$.

If in addition, $G$ is a RIPL or LIPL then by earlier result and
since $J_\rho^2=I$ or $J_\lambda^2=I$ respectively, $C=D$ in each
case. Consequently, $J_\rho =J_\lambda$.

\begin{myrem}
In \cite{phd30}, isotopic CIP loops are shown to be isomorphic and
this is true for commutative Moufang loops as shown in \cite{phd3}.
\end{myrem}

\section{Conclusion and Future Study}
Karklin$\ddot{\textrm{u}}$sh and Karkli$\check{\textrm{n}}$
\cite{phd175} introduced $m$-inverse loops i.e loops that obey any
of the equivalent conditions
\begin{displaymath}
(xy)J_\rho^m\cdot xJ_\rho^{m+1}=yJ_\rho^m\qquad\textrm{and}\qquad
xJ_\lambda^{m+1}\cdot (yx)J_\lambda^m=yJ_\lambda^m.
\end{displaymath}
They are generalizations of WIPLs and CIPLs, which corresponds to
$m=-1$ and $m=0$ respectively. After the study of $m$-loops by
Keedwell and Shcherbacov \cite{phd176}, they have also generalized
them to quasigroups called $(r,s,t)$-inverse quasigroups in
\cite{phd177} and \cite{phd178}. It will be interesting to study the
universality of $m$-inverse loops and $(r,s,t)$-inverse quasigroups.
These will generalize the works of J. M. Osborn and R. Artzy on
universal WIPLs and CIPLs respectively. Furthermore, we raise the
question of studying $m$-inverse loops and $(r,s,t)$-inverse
quasigroups under a generalized type of ${\cal T}$ condition say a
${\cal T}_m$ condition and ${\cal T}_{(r,s,t)}$ condition just like
we have done it for WIPLs? Also, we need to know if under the ${\cal
T}$ condition, a loop in a pair is a CIPL if and only if the other
loop is a CIPL as we have shown it to be true for WIPLs in
Theorem~\ref{1:3}.

\end{document}